\documentclass[11pt,a4paper]{article}
\usepackage[T1]{fontenc}
\usepackage[utf8]{inputenc}
\usepackage{lmodern}
\usepackage[margin=25mm]{geometry}
\usepackage{amsmath,amssymb,amsthm}
\usepackage{graphicx,booktabs,multirow}
\usepackage{subfig}
\usepackage{float}
\usepackage[linesnumbered,ruled,vlined]{algorithm2e}
\usepackage[numbers,sort&compress]{natbib}

\usepackage{microtype}
\usepackage[hidelinks]{hyperref}
\hypersetup{
  pdftitle={Efficient numerical computation of the generalized Marcum Q-function with applications},
  pdfauthor={Zhong-Xuan Mao and Jing-Feng Tian}
}

\theoremstyle{plain}
\newtheorem{theorem}{Theorem}[section]
\newtheorem{corollary}[theorem]{Corollary}
\theoremstyle{definition}
\newtheorem{remark}{Remark}[section]

\title{Efficient numerical computation of the generalized Marcum $Q$-function with applications}
\author{Zhong-Xuan Mao\textsuperscript{1,2}
\and Jing-Feng Tian\textsuperscript{2}\thanks{Corresponding author. Email: \texttt{tianjf@ncepu.edu.cn}.}}
\date{}

\begin{document}
\maketitle
\vspace{-2em}
\begin{center}
\small
\textsuperscript{1}Jilin University, Changchun, China\\
\textsuperscript{2}North China Electric Power University, Baoding, China
\end{center}

\begin{abstract}
The computation of the generalized Marcum $Q$-function, $Q_\nu(a,b)$, is fundamental to radar detection and wireless communications. The state-of-the-art \textbf{Algorithm 939} (Gil et al., ACM TOMS 2014) serves as the current gold standard by employing a hybrid evaluation scheme combining series expansions, asymptotic formulas,  numerical integration and recurrence relations. However, it suffers from a computational bottleneck: in the ``large $a$, small $b$'' regime, the algorithm resorts to computationally expensive numerical quadrature, thereby limiting its overall efficiency.

To bridge this gap, we utilize a series representation of the generalized Marcum $Q$-function to construct the Series Expansion based on Recurrence Relations (\textbf{SE-RR}) algorithm. Expanded in powers of $b$, this series computes its coefficients via a stable recurrence relation, resulting in a remarkably low computational cost. By analyzing rounding, truncation and cancellation errors, we identify the optimal region for SE-RR as the ``large $a$, small $b$'' regime, which covers the parameter space where Algorithm~939 relies on costly numerical quadrature. Consequently, we integrate SE-RR into Algorithm~939 to develop a hybrid algorithm for evaluating $Q_\nu(a,b)$. Benchmarked against $10^7$ log-uniformly sampled parameter triples $(a, b, \nu) \in [0.01, 100] \times [0.001, 100] \times [1, 1000]$, the hybrid algorithm achieves an average speedup of 82.2\% overall, while in the region where SE-RR is applied, it delivers an average speedup of 156.6\% and a maximum pointwise speedup of $\mathrm{5.49}\times$.
Finally, we demonstrate the practical utility of the hybrid algorithm by computing the cumulative distribution function of the noncentral chi distribution and outage probabilities for 5G/6G ultra-reliable low-latency communications.
\end{abstract}
\noindent\textbf{Keywords:} Generalized Marcum $Q$-function, Noncentral chi distribution, Modified Bessel function of the first kind, Series expansion
\par
\noindent\textbf{2020 Mathematics Subject Classification:} Primary 65D20; Secondary 33C10, 33C45, 65G50.\par

\section{Introduction}

The generalized Marcum $Q$-function of order $\nu>0$ is defined by
\begin{equation} \label{Qab}
Q_\nu(a,b) = \frac{1}{a^{\nu-1}} \int_{b}^\infty t^\nu e^{-\frac{t^2+a^2}{2}}I_{\nu-1}(at)  \textrm{d}t, \quad a>0,\; b\geq 0,
\end{equation}
where $I_\nu$ is the modified Bessel function of the first kind defined by
\begin{equation} \label{I-series}
I_\nu(x) = \sum_{n=0}^\infty \frac{1}{n! \Gamma(\nu+n+1)} \Big(\frac{x}{2}\Big)^{2n+\nu}, \quad \nu \notin \{-1,-2,\cdots\}.
\end{equation}
In particular, when $\nu=1$, \eqref{Qab} reduces to the well-known Marcum $Q$-function
\begin{equation*}
Q(a,b) := Q_1(a,b) = \int_b^\infty t e^{-\frac{t^2+a^2}{2}} I_{0}(at) \textrm{d}t.
\end{equation*}
The generalized Marcum $Q$-function is closely related to the noncentral chi distribution. Specifically, if $R$ follows a noncentral chi distribution with $2M$ degrees of freedom and noncentrality parameter $a$, then $\mathbb{P}(R>b) = Q_M(a,b)$.
Moreover, $Q_\nu(a,b)$ is widely used in radar and wireless communications, with key applications including pulse radar target detection \cite{Marcum-IRE-TIT-1960}, error rate analysis over fading channels \cite{Simon-IEEE-TC-2000}, performance evaluation of predictor antenna systems \cite{Guo-IEEE-OJCS-2021}, and average detection probability analysis for energy detectors \cite{Sofotasios-IEEE-SPL-2015,Cao-IEEE-CL-2016}.

Consequently, the efficient and accurate computation of this function is of significant importance. One category of methods utilizes the analytic properties of the modified Bessel function to establish lower and upper bounds  on \eqref{Qab}, as seen in the works of Baricz and Sun \cite{Baricz-IEEE-TIT-2009}, Baricz \cite{Baricz-JMAA-2009}, and others \cite{Simon-IEEE-TC-2000,Corazza-IEEE-TIT-2002,Baricz-JCAM-2021,Baricz-RJ-2022,Mao-AMC-2027}. This methodology is straightforward and intuitive, inherently providing an estimate for the maximum error, which is defined by the difference between the bounds. While this approach offers an intuitive error estimate, achieving high accuracy in practical computations often requires further refinement.

Other researchers have addressed the problem from a computational perspective. The most direct strategy involves the numerical integration of \eqref{Qab}, for which numerous techniques are available. In particular, Temme \cite{Temme-CMA-1993} derived a contour integral representation for $Q_\nu(a,b)$, based on which Gil et al. \cite{Gil-ACM-2014} developed a comprehensive algorithm. Although numerical integration algorithms are the most robust---capable of handling various extreme parameter regions---they suffer from slower computation speeds, as techniques such as the trapezoidal rule require multiple evaluations of the integrand.
Alternatively, the series expansion method utilizes the series expansion of the modified Bessel function~\eqref{I-series} to yield
\begin{equation} \label{QS}
Q_\nu(a,b) = e^{-\frac{a^2}{2}} \sum_{n=0}^\infty \frac{1}{n!} \frac{\Gamma\big(\nu+n, \frac{b^2}{2}\big)}{\Gamma(\nu+n)} \Big(\frac{a^2}{2}\Big)^n ,
\end{equation}
where $\Gamma(s,x)= \int_x^\infty u^{s-1} e^{-u} \textrm{d} u$ denotes the upper incomplete gamma function.

This method boasts the fastest computation speed among existing techniques, but it is limited to handling relatively small values of $a$. Addressing the case of large $a$ and $b$, Gil et al. \cite{Gil-ACM-2014} established two sets of asymptotic expansion algorithms designed for small and large $\nu$, respectively. These asymptotic methods offer moderate computation speeds and are highly efficient for large parameters. Furthermore, Temme \cite{Temme-CMA-1993} derived recurrence relations for $Q_\nu(a,b)$ with respect to $\nu$, which are utilized for numerical computation in specific transition regions. This algorithm essentially maps points from one region to another, and in practice, it still relies on alternative methods to provide an initial value to seed the recurrence.

Based on the characteristics of the aforementioned methods, Gil et al.~\cite{Gil-ACM-2014} proposed a hybrid algorithm (known as \textbf{Algorithm 939}) for the efficient computation of $Q_\nu(a,b)$.
When $\nu$ is not extremely large, three primary computational methods are employed: series expansions (\textbf{SE}) are used for small values of $a$; numerical integration (\textbf{NI}) is applied when $b$ is small and $a$ is large; and asymptotic formulas (\textbf{AF}) are adopted for large values of both $a$ and $b$.
Near the curve $b^2=a^2+2\nu$, Algorithm 939 also use recurrence relations (\textbf{RR}) for $\nu<135$ and a large-order asymptotic expansion (\textbf{LAE}) for large $\nu$.

However, a critical gap remains in the existing computational framework. In the regime of small $b$ and large $a$, one is forced to rely on numerical integration, whose computational efficiency is typically much lower than that of series-based evaluation. This lack of an efficient series-type method in this parameter region constitutes a significant bottleneck in practical computations.

To address this limitation, we utilize a series representation of the generalized Marcum $Q$-function to construct the \textbf{SE-RR} algorithm. By expanding in powers of $b$ and computing coefficients via a stable recurrence relation, SE-RR achieves high computational efficiency. Guided by an analysis of rounding error, cancellation, and truncation error, we demonstrate that SE-RR effectively targets the ``large $a$, small $b$'' regime, where existing methods rely on costly numerical quadrature. Incorporating SE-RR into Algorithm~939 yields an improved hybrid algorithm. Benchmarked against $10^7$ log-uniformly sampled parameter triples $(a,b,\nu)\in[0.01,100]\times[0.001,100]\times[1,1000]$, this hybrid approach achieves an average speedup of $82.2\%$ overall and $156.6\%$ in the region where SE-RR is applied, with a maximum pointwise speedup of $5.49\times$.

Finally, we demonstrate the practical value of the proposed SE-RR algorithm in the chi distribution and in 5G/6G Ultra-Reliable Low-Latency Communications (URLLC). We present an efficient algorithm for computing the Cumulative Distribution Function (CDF) $F_R(b_i)$ $(i=1,\cdots,n)$  of $R$, where $R$ follows a noncentral chi distribution with $2M$ degrees of freedom and noncentrality parameter $a$. Since the coefficients in the proposed SE-RR are independent of $b$, they only need to be computed once rather than being recomputed at each node $b_i$.
In high-reliability URLLC scenarios, the wireless channel typically exhibits a strong Line-of-Sight (LoS) component, corresponding to a large Rician $K$-factor. On the other hand, to guarantee ultra-high reliability, the threshold $\gamma_{\text{th}}$ is usually chosen to be very small. Consequently, for a normalized average SNR $\bar{\gamma}=1$, the outage probability
\begin{equation*}
P_{\text{out}} = 1 - Q_1\big( \sqrt{2K}, \sqrt{2(1+K)\gamma_{\text{th}}} \big)
\end{equation*}
requires the evaluation of a generalized Marcum $Q$-function with a relatively large $a$ and a small $b$, a regime in which the proposed algorithm is particularly advantageous.

\section{Preliminaries}
In this section, we first review several fundamental properties of the generalized Marcum $Q$-function. We then briefly introduce SE, AF, and NI, and explain how Algorithm~939 selects among these methods alongside RR and LAE across different parameter regions.

Clearly, the functions $a \mapsto Q_\nu(a,b)$ and $\nu \mapsto Q_\nu(a,b)$ are increasing, whereas $b \mapsto Q_\nu(a,b)$ is decreasing~\cite{Gil-ACM-2014}.
Noting that $Q_\nu(a,0)=1$~\cite{Baricz-IEEE-TIT-2009} and $\lim_{b\to\infty} Q_\nu(a,b)=0$, we obtain
$0 \leq Q_\nu(a,b) \leq 1$ for all $a>0$, $b\geq0$, and $\nu>0$.
Moreover, for all $\nu>0$, the complementary integral representation~\cite{Baricz-JMAA-2009} is
\begin{equation} \label{Qab-symmetry}
Q_\nu(a,b) = 1 - \frac{1}{a^{\nu-1}} \int_{0}^b t^\nu e^{-\frac{t^2+a^2}{2}} I_{\nu-1}(at) \, dt,
\quad a>0,\; b\geq 0.
\end{equation}
In some works, for example \cite{Gil-ACM-2014}, the complementary function $P_\nu(a,b):=1 - Q_\nu(a,b)$ is also studied extensively. In particular, when $Q_\nu(a,b)$ is close to $1$, evaluating $1-Q_\nu(a,b)$ directly may suffer from severe cancellation errors, and therefore it is often numerically preferable to compute $P_\nu(a,b)$ directly instead.

Next, we introduce the SE, AF, and NI methods employed in Algorithm 939.
We have already derived the SE representation in \eqref{QS}.
For the coefficients $\Gamma(\nu+n, \frac{b^2}{2})/\Gamma(\nu+n)$, the incomplete gamma function does not need to be evaluated repeatedly.
Instead, one may employ the recurrence relations $\Gamma(a+1,x) = a \Gamma(a,x) + x^a e^{-x}$ and $\Gamma(\nu+n+1)= (\nu+n)\Gamma(\nu+n)$,
so that only the initial value $\Gamma(\nu, \frac{b^2}{2})$ and $\Gamma(\nu)$ needs to be computed explicitly.
Alternatively, one may use the specialized algorithm proposed in \cite{Gil-SIAM-JSC-2012}, as recommended in \cite{Gil-ACM-2014}, for the efficient evaluation of functions of the form
$\Gamma(a,x)/\Gamma(a)$.

By exploiting the asymptotic expansion of the modified Bessel function, we obtain the AF of the generalized Marcum $Q$-function.
For large $ab$, it is given by
\begin{equation*}
Q_\nu(a,b) \sim \sum_{n=0}^\infty \frac{(b/a)^\nu}{2\sqrt{2\pi}} (-1)^n \big( A_n(\nu - 1) - \frac{a}{b} A_n(\nu) \big) \Phi_n,
\end{equation*}
where
\begin{equation*}
A_n(\nu) = \frac{\Gamma(\frac{1}{2}+\nu+n)}{2^n n! \Gamma(\frac{1}{2}+\nu-n)},
\end{equation*}
and
\begin{equation*}
\Phi_n = \Big( \frac{(a-b)^2}{2ab} \Big)^{n-\frac{1}{2}} \Gamma\Big(\frac{1}{2}-n,\frac{(a-b)^2}{2} \Big).
\end{equation*}
The above method is applicable for the case $a < b$. When $a \geq b$, the complementary function is computed by
\begin{equation*}
P_\nu(a,b) \sim \sum_{n=0}^\infty \frac{(b/a)^\nu}{2\sqrt{2\pi}} (-1)^{n+1} \big( A_n(\nu - 1) - \frac{a}{b} A_n(\nu) \big) \Phi_n,
\end{equation*}

Next, we introduce NI, which is given by
\begin{equation*}
Q_\nu(a,b)=
\begin{cases}
\displaystyle \frac{e^{-\frac{1}{2}\nu\zeta^2}}{2\pi} \int_{-\pi}^{\pi} e^{\nu\psi(\theta)} f(\theta) \, \mathrm{d}\theta,
& b^2>a^2+2\nu,\\[8pt]
\displaystyle 1+\frac{e^{-\frac{1}{2}\nu\zeta^2}}{2\pi} \int_{-\pi}^{\pi} e^{\nu\psi(\theta)} f(\theta) \, \mathrm{d}\theta,
& b^2<a^2+2\nu,
\end{cases}
\end{equation*}
where $\xi = ab/\nu$,
\begin{equation*}
\zeta^2 = 2 \bigl( \frac{a^2+b^2}{2\nu} - \sqrt{1+\xi^2} + \ln\bigl( \frac{1+\sqrt{1+\xi^2}}{b^2/\nu} \bigr) \bigr),
\end{equation*}
and its sign is consistent with $a^2+2\nu - b^2$.
The functions $\rho$ and $r$ are defined by
\begin{equation*}
\rho(\theta) = \sqrt{\bigl(\frac{\theta}{\sin\theta}\bigr)^2 + \xi^2}, \quad
r(\theta) = \frac{\nu}{b^2} \bigl( \frac{\theta}{\sin\theta} + \rho(\theta) \bigr),
\end{equation*}
and $f$ and $\psi$ are given by
\begin{equation*}
f(\theta) = \frac{\sin\theta \, r'(\theta) + (\cos\theta - r(\theta)) r(\theta)}{r^2(\theta) - 2 r(\theta)\cos\theta + 1}, \quad
\psi(\theta) = \rho(\theta)\cos\theta - \sqrt{1+\xi^2} - \ln\bigl( \frac{\frac{\theta}{\sin\theta} + \rho(\theta)}{1+\sqrt{1+\xi^2}} \bigr).
\end{equation*}

\begin{figure}[htbp]
    \centering
    \subfloat[$\nu=1$]{
        \includegraphics[width=0.45\textwidth]{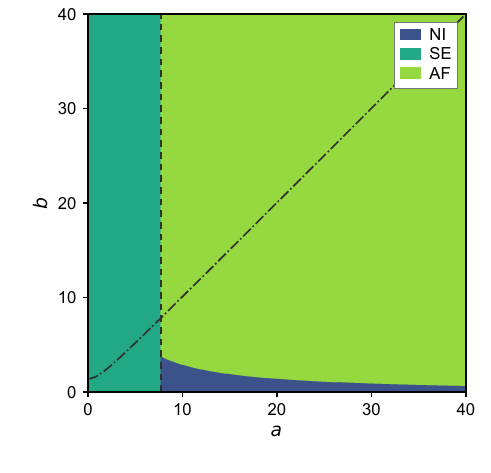}
        \label{fig:nu1}
    }
    \hfill
    \subfloat[$\nu=20$]{
        \includegraphics[width=0.45\textwidth]{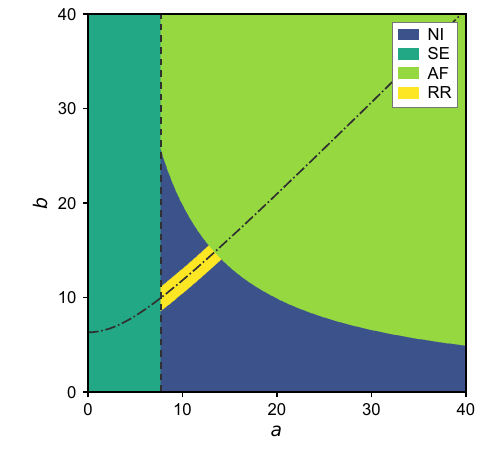}
        \label{fig:nu10}
    }
    \caption{The regions in the positive quadrant of the $(a,b)$-plane where specific algorithms are applied for $\nu=1$ and $\nu=20$, based on Algorithm 939.}
    \label{fig:ab_regions}
\end{figure}

The applicability and computational efficiency of the above three methods differ significantly across parameter regimes.
To address this issue, Gil et al.~\cite{Gil-ACM-2014} proposed a hybrid strategy for the efficient evaluation of $Q_\nu(a,b)$.
The strategy also uses RR, which computes the Marcum functions via a three-term recurrence in the order $\nu$, together with LAE, which approximates them by applying a large-order expansion of the modified Bessel function to their integral representations.
Define
\begin{equation*}
    g_{1,2}(a, \nu) = \sqrt{a^2 + 2\nu \mp 2\sqrt{2a^2 + 2\nu}}.
\end{equation*}
Algorithm 939 is stated as follows:
\begin{enumerate}
    \item If $a < \sqrt{60} \approx 7.746$, use SE.
    \item If $ab > 30$ and $\nu^2 < 2ab$, use AF.
    \item If $g_1(a, \nu) < b < g_2(a, \nu)$ and $\nu < 135$, use RR.
    \item If $g_1(a, \nu) < b < g_2(a, \nu)$ and $\nu \ge 135$, use LAE.
    \item Otherwise, use NI.
\end{enumerate}
When $\nu$ is not extremely large, three primary computational methods are employed (see Figure \ref{fig:ab_regions}): series expansions are used for small values of $a$; numerical integration is applied when $b$ is small and $a$ is large; and asymptotic formulas are adopted for large values of both $a$ and $b$.

\section{Recurrence relations for the expansion coefficients of the product of an exponential function and a modified Bessel function}
\label{sec3}

We first outline the main idea of SE-RR. Suppose that $e^{-t^2/2} I_{\nu-1}(at)$ admits the expansion
\begin{equation*}
e^{-\frac{t^2}{2}} I_{\nu-1}(at) = \sum_{n=0}^\infty v_n(a,\nu) t^{2n+\nu-1},
\end{equation*}
where the coefficients $v_n(a,\nu)$ can be computed via a recurrence relation. Substituting this expansion into \eqref{Qab-symmetry} yields
\begin{equation*}
Q_\nu(a,b) = 1 - \frac{e^{-\frac{a^2}{2}}}{a^{\nu-1}} \sum_{n=0}^\infty v_n(a,\nu) \frac{b^{2n+2\nu}}{2n+2\nu}.
\end{equation*}
Up to the prefactor $b^{2\nu}$, the sum is a power series in $b^2$, making this representation particularly effective for small values of $b$. Thus, we first consider the series expansion of
\begin{equation*}
t \mapsto \exp\!\big(-t^2/p\big) I_{\mu}(at),
\end{equation*}
where $p,a,\mu \in \mathbb{R}$ and $\mu \notin \{0,-1,-2,\ldots\}$.

\begin{theorem} \label{thm-main}
Let $a, \mu, p ,t \in \mathbb{R}$ with $p\neq 0$ and $\mu \notin \{-1, -2, \dots\}$. Then
\begin{equation*}
\exp\big(-\frac{t^2}{p}\big) I_{\mu}(at) = \frac{a^{\mu}}{2^{\mu} \Gamma(\mu+1)} \sum_{n=0}^\infty u_n(a,\mu,p) t^{2n+\mu},
\end{equation*}
where the coefficients are determined by the initial values $u_0=1$ and $u_1=\frac{a^2}{4(\mu+1)}-\frac{1}{p}$ together with the following three-term recurrence relation
\begin{equation*}
u_n = \frac{a^2 p - 4(\mu + 2n - 1)}{4pn(\mu+n)} u_{n-1} - \frac{1}{p^2 n(\mu+n)} u_{n-2}, \quad n\geq 2.
\end{equation*}
\end{theorem}

\begin{proof}
Let $C_{a,\mu} := \frac{a^{\mu}}{2^{\mu} \Gamma(\mu+1)}$ and define
\begin{equation*}
F(t) = \exp\big(-\frac{t^2}{p}\big) I_{\mu}(at).
\end{equation*}
Suppose that
\begin{equation*}
F(t) = C_{a,\mu} \sum_{n=0}^\infty u_n t^{2n+\mu}.
\end{equation*}
To establish the recurrence relation for the coefficients $u_n$, it suffices to verify that the following equation holds for all $t$
\begin{equation} \label{ALL}
C_{a,\mu} \sum_{n=0}^\infty \Big( n(\mu+n)u_n - \frac{a^2 p - 4(\mu + 2n - 1)}{4p} u_{n-1} + \frac{1}{p^2} u_{n-2} \Big) t^{2n+\mu} = 0,
\end{equation}
where we assume $u_{-1} = u_{-2} = 0$.

We analyze the terms in the summation individually by relating them to the derivatives of $F(t)$. First, observe that
\begin{equation*}
\begin{aligned}
C_{a,\mu} \sum_{n=0}^\infty n(\mu+n) u_n t^{2n+\mu} & = C_{a,\mu} \sum_{n=0}^\infty \Big( \frac{(2n+\mu)(2n+\mu-1)}{4} + \frac{2n+\mu}{4} - \frac{\mu^2}{4} \Big) u_n t^{2n+\mu} \\
& = \frac{t^2}{4} F''(t) + \frac{t}{4} F'(t) - \frac{\mu^2}{4} F(t).
\end{aligned}
\end{equation*}
Next, for the second term, we have
\begin{equation*}
\begin{aligned}
- C_{a,\mu} \sum_{n=0}^\infty \frac{a^2 p - 4(\mu + 2n - 1)}{4p} u_{n-1} t^{2n+\mu} & = C_{a,\mu} \sum_{n=0}^\infty \Big( \frac{2(n-1)+\mu}{p} + \frac{1}{p}-\frac{a^2}{4} \Big) u_{n-1} t^{2n+\mu} \\
& = \frac{t^3}{p} F'(t) + \Big( \frac{1}{p}-\frac{a^2}{4} \Big) t^2 F(t).
\end{aligned}
\end{equation*}
Finally, the third term corresponds to
\begin{equation*}
C_{a,\mu} \sum_{n=0}^\infty \frac{1}{p^2} u_{n-2} t^{2n+\mu} = \frac{t^4}{p^2} F(t).
\end{equation*}
Substituting these relations into the left-hand side of \eqref{ALL}, we obtain the differential expression
\begin{equation} \label{LHS}
\frac{t^2}{4} F^{\prime\prime}(t) + \Big( \frac{t}{4} + \frac{t^3}{p} \Big) F^{\prime}(t) + \Big(- \frac{\mu^2}{4} + \big( \frac{1}{p}-\frac{a^2}{4} \big) t^2 + \frac{t^4}{p^2}\Big)F(t).
\end{equation}
To verify that \eqref{LHS} equals zero, we compute the derivatives of $F(t)$. Using the property of the modified Bessel function
\begin{equation*}
I_{\mu}^{\prime\prime}(at) = -\frac{1}{t} I_{\mu}^{\prime}(at) + \Big(a^2+\frac{\mu^2}{t^2}\Big) I_{\mu}(at),
\end{equation*}
the first derivative is
\begin{equation} \label{F_prime}
F^{\prime}(t) = -\frac{2t}{p} e^{-\frac{t^2}{p}} I_{\mu}(at) + e^{-\frac{t^2}{p}} I_{\mu}^{\prime}(at).
\end{equation}
Differentiating again, we find
\begin{equation} \label{F_double_prime}
\begin{aligned}
F^{\prime\prime}(t) &= \frac{\mathrm{d}}{\mathrm{d} t} \Big( -\frac{2t}{p} F(t) + e^{-\frac{t^2}{p}} I_{\mu}'(at) \Big) \\
&= \Big(\frac{4 t^2}{p^2}-\frac{2}{p}\Big) e^{-\frac{t^2}{p}} I_{\mu}(at) - \frac{4 t}{p} e^{-\frac{t^2}{p}} I_{\mu}^{\prime}(at) + e^{-\frac{t^2}{p}} I_{\mu}^{\prime\prime}(at) \\
&= \frac{e^{-\frac{t^2}{p}}}{p^2 t^2} ( I_{\mu}(at)( p t^2 (a^2 p - 2) + \mu^2 p^2 + 4 t^4 ) - p t (p + 4 t^2) I_{\mu}^{\prime}(at) ).
\end{aligned}
\end{equation}
Substituting \eqref{F_prime} and \eqref{F_double_prime} into \eqref{LHS}, algebraic simplification confirms that the expression vanishes identically. Thus, the recurrence relation holds.

The initial coefficients $u_0$ and $u_1$ can be determined directly using the Cauchy product of the series expansions for $\exp(-t^2/p)$ and $I_\mu(at)$.
\end{proof}

\begin{remark}
\label{rem:laguerre-coefficients}
As pointed out by an anonymous reviewer, the coefficients in Theorem~\ref{thm-main} admit an explicit representation in terms of generalized Laguerre polynomials
\begin{equation}
\label{un-Laguerre}
u_n(a,\mu,p) = \frac{(-1/p)^n}{(\mu+1)_n} L_n^{(\mu)}\Bigl(\frac{a^2p}{4}\Bigr),
\end{equation}
where $(\mu+1)_n$ denotes the Pochhammer symbol and $L_n^{(\mu)}$ is the generalized Laguerre polynomial of degree $n$, defined by \cite[Eq.~(18.5.12)]{Olver-2010}
\begin{equation*}
L_n^{(\mu)}(x) = \sum_{k=0}^{n} \frac{(-1)^k\Gamma(n+\mu+1)}{\Gamma(\mu+k+1)(n-k)!\,k!} x^k.
\end{equation*}
Consequently, the recurrence relation in Theorem~\ref{thm-main} can also be established directly from the standard three-term recurrence for Laguerre polynomials \cite[Eq.~(18.9.1)]{Olver-2010}
\begin{equation*}
n L_n^{(\mu)}(x) = (2n+\mu-1-x) L_{n-1}^{(\mu)}(x) - (n+\mu-1) L_{n-2}^{(\mu)}(x), \qquad n \ge 2.
\end{equation*}
\end{remark}

Taking $p=2$ in Theorem \ref{thm-main}, we obtain the following corollary.

\begin{corollary} \label{cor-main}
Let $a, \mu, t \in \mathbb{R}$ with $\mu \notin \{-1, -2, \dots\}$. Then
\begin{equation*}
\exp\big(-\frac{t^2}{2}\big) I_{\mu}(at) = \frac{a^{\mu}}{2^{\mu} \Gamma(\mu+1)} \sum_{n=0}^\infty u_n(a,\mu) t^{2n+\mu},
\end{equation*}
where $u_0=1$, $u_1=\frac{a^2 - 2(\mu +1)}{4(\mu +1)}$ and
\begin{equation*}
u_n = \frac{a^2 - 2(\mu + 2n - 1)}{4n(\mu+n)} u_{n-1} - \frac{1}{4n(\mu+n)} u_{n-2}, \quad n\geq 2.
\end{equation*}
\end{corollary}

Taking $\mu=0$ in Corollary \ref{cor-main}, we obtain the following corollary.
\begin{corollary} \label{cor-0}
Let $a, t \in \mathbb{R}$. Then
\begin{equation*}
\exp\big(-\frac{t^2}{2}\big) I_{0}(at) = \sum_{n=0}^\infty u_n(a) t^{2n},
\end{equation*}
where $u_0=1$, $u_1=\frac{a^2 - 2}{4}$, and the coefficients satisfy the recurrence relation
\begin{equation*}
u_n = \frac{a^2 - 2(2n - 1)}{4n^2} u_{n-1} - \frac{1}{4n^2} u_{n-2}, \quad n\geq 2.
\end{equation*}
\end{corollary}

\begin{figure}[htbp]
    \centering
        \includegraphics[width=0.8\textwidth]{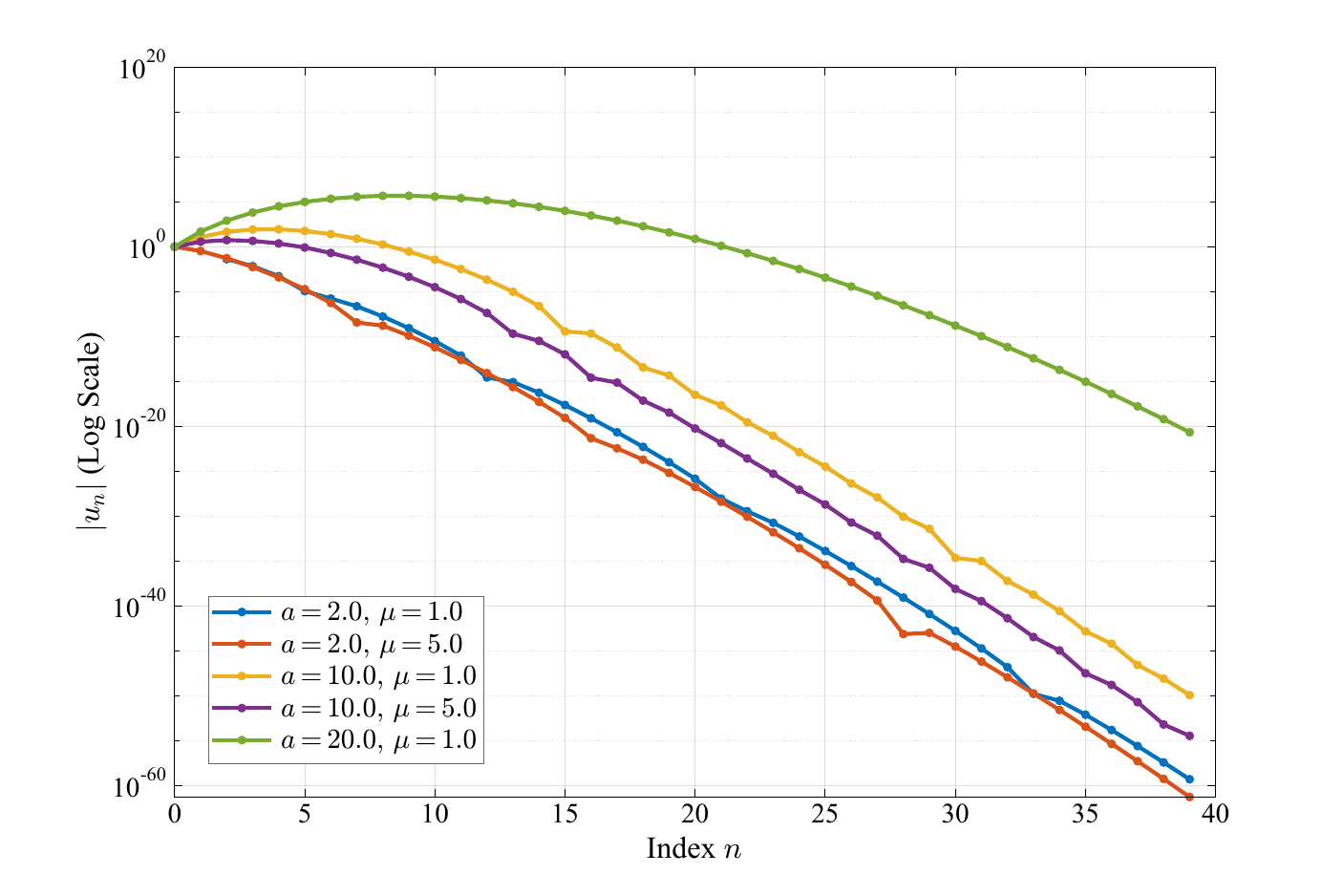}
    \caption{Decay trend of coefficients $|u_n|$}
    \label{fig-un}
\end{figure}

Next, we investigate the properties of the coefficients determined by the recurrence relation in Corollary~\ref{cor-main}. To begin with, we select several representative parameter settings and illustrate the behavior of $u_n(a,\nu)$ through numerical plots, see Figure \ref{fig-un}.

In fact, we can establish the following asymptotic behavior of $u_n$, which indicates that $u_n \to 0$ as $n \to \infty$ at an exponential rate.

\begin{theorem} \label{thm-ab}
Let $a>0$ and $\nu \notin \{0, -1, -2, \dots\}$. Suppose that $u_0=1$, $u_1=\frac{a^2 - 2\nu}{4\nu}$ and
\begin{equation}\label{eq:un-recurrence}
u_n = \frac{a^2 - 2(\nu + 2n - 2)}{4n(\nu+n-1)} u_{n-1} - \frac{1}{4n(\nu+n-1)} u_{n-2}, \quad n\geq 2.
\end{equation}

Then
\begin{equation*}
u_n = \mathcal{O}\Bigl( \frac{n^{\frac{1}{4}-\frac{\nu}{2}}}{2^n n!} \Bigr), \quad \text{as } \, n\to \infty.
\end{equation*}

\end{theorem}

\begin{proof}
To facilitate the asymptotic analysis, we introduce the auxiliary sequence $w_n = n! u_n$. Substituting this transformation into the recurrence relation yields
\begin{equation} \label{w-rr}
w_{n+1} = \frac{a^2 - 2(\nu + 2n)}{4(\nu+n)} w_n - \frac{n}{4(\nu+n)} w_{n-1}, \qquad n\ge1.
\end{equation}
As $n \to \infty$, a straightforward calculation of the coefficients reveals their limiting values
\begin{equation*}
\lim_{n \to \infty} \frac{a^2 - 2(\nu + 2n)}{4(\nu+n)} = -1 \quad \text{and} \quad \lim_{n \to \infty} - \frac{n}{4(\nu+n)} = -\frac{1}{4}.
\end{equation*}
This leads to the characteristic equation $t^2 + t + \frac{1}{4}=0$, which possesses a repeated root of multiplicity two at $\lambda_{1,2} = -\frac{1}{2}$.

The Birkhoff--Trjitzinsky theory for difference equations \cite[Theorem~1]{Wimp-JMAA-1985} suggests that $w_n$ possesses the asymptotic form
\begin{equation*}
w_n\sim e^{\mu_0n\log n+\mu_1n}n^\theta\exp\big(\alpha_1n^\beta+\alpha_2n^{\beta-1/\rho}+\alpha_3n^{\beta-2/\rho}+\cdots\big),
\end{equation*}
where $\alpha_1\neq0$, $0<\beta<1$, $\rho,\beta\rho\in\mathbb Z_{>0}$, and $e^{\mu_0+\mu_1}=\lambda=-\frac12$.
Balancing the dominant terms associated with the factor $e^{\mu_0 n \log n + \mu_1 n}$ in the recurrence relation \eqref{w-rr} yields $\mu_0 = 0$ and $e^{\mu_1} = -1/2$. Consequently, setting $\Phi(n) = \alpha_1 n^\beta + \alpha_2 n^{\beta-1/\rho} + \alpha_3 n^{\beta-2/\rho} + \cdots$, we obtain
\begin{equation*}
w_n \sim \Bigl(-\frac{1}{2}\Bigr)^n n^\theta e^{\Phi(n)}.
\end{equation*}
Then, formally, the shift ratios are given by
\begin{equation*}
\frac{w_{n+2}}{w_n} = \frac{1}{4} \Bigl(1+\frac{2}{n}\Bigr)^\theta e^{\Phi(n+2)-\Phi(n)}, \qquad
\frac{w_{n+1}}{w_n} = -\frac{1}{2} \Bigl(1+\frac{1}{n}\Bigr)^\theta e^{\Phi(n+1)-\Phi(n)}.
\end{equation*}
Substituting these ratios into the recurrence relation \eqref{w-rr} gives
\begin{equation}
\begin{aligned}
0 ={}& (n+\nu+1) \Bigl(1+\frac{2}{n}\Bigr)^\theta e^{\Phi(n+2)-\Phi(n)} \\
& - \Bigl(2n+\nu+2-\frac{a^2}{2}\Bigr) \Bigl(1+\frac{1}{n}\Bigr)^\theta e^{\Phi(n+1)-\Phi(n)} + n + 1.
\end{aligned}
\label{eq:BT-substitution}
\end{equation}
Moreover, we have the asymptotic expansions
\begin{equation*}
\begin{aligned}
\big(1+\frac{2}{n}\big)^\theta &= 1 + \frac{2\theta}{n} + \mathcal{O}(n^{-2}), \\
\big(1+\frac{1}{n}\big)^\theta &= 1 + \frac{\theta}{n} + \mathcal{O}(n^{-2}), \\
e^{\Phi(n+2)-\Phi(n)} &= 1 + 2\alpha_1\beta n^{\beta-1} + 2\alpha_1^2\beta^2 n^{2\beta-2} + \mathcal{O}\bigl(n^{\beta-1/\rho-1} + n^{\beta-2} + n^{3\beta-3}\bigr), \\
e^{\Phi(n+1)-\Phi(n)} &= 1 + \alpha_1\beta n^{\beta-1} + \frac{1}{2}\alpha_1^2\beta^2 n^{2\beta-2} + \mathcal{O}\bigl(n^{\beta-1/\rho-1} + n^{\beta-2} + n^{3\beta-3}\bigr).
\end{aligned}
\end{equation*}
Substituting these expansions into \eqref{eq:BT-substitution} yields the leading-order balance
\begin{equation*}
0 = \alpha_1^2\beta^2 n^{2\beta-1} + \frac{a^2}{2} + o\bigl(n^{2\beta-1}\bigr) + o(1),
\end{equation*}
which implies $\beta = 1/2$ and $\alpha_1 = \pm\mathrm{i}\sqrt{2}\,a$. By setting $\rho = 2$, we obtain the refined expansions
\begin{equation*}
\begin{aligned}
e^{\Phi(n+2)-\Phi(n)} &= 1 + \alpha_1 n^{-1/2} + \frac{\alpha_1^2}{2}n^{-1} + \Bigl(\frac{\alpha_1^3}{6} - \frac{\alpha_1}{2} - \alpha_3\Bigr)n^{-3/2} + \mathcal{O}(n^{-2}), \\
e^{\Phi(n+1)-\Phi(n)} &= 1 + \frac{\alpha_1}{2}n^{-1/2} + \frac{\alpha_1^2}{8}n^{-1} + \Bigl(\frac{\alpha_1^3}{48} - \frac{\alpha_1}{8} - \frac{\alpha_3}{2}\Bigr)n^{-3/2} + \mathcal{O}(n^{-2}).
\end{aligned}
\end{equation*}
Substituting these refined expansions back into \eqref{eq:BT-substitution} again gives
\begin{equation*}
0 = \frac{\alpha_1^2}{4} + \frac{a^2}{2} + \frac{\alpha_1}{8} \bigl(\alpha_1^2+2a^2+4\nu+8\theta-2\bigr) n^{-1/2} + \mathcal{O}(n^{-1}),
\end{equation*}
which implies $\theta = \frac{1}{4} - \frac{\nu}{2}$. Consequently, the general asymptotic solution takes the form
\begin{equation*}
\begin{aligned}
w_n = \Bigl(-\frac{1}{2}\Bigr)^n n^{\frac{1}{4}-\frac{\nu}{2}} \Bigl(
&C_{+} \exp\bigl(\mathrm{i}\sqrt{2}\,a\sqrt{n}\bigr) \bigl(1+\mathcal{O}(n^{-1/2})\bigr) \\
+{}& C_{-} \exp\bigl(-\mathrm{i}\sqrt{2}\,a\sqrt{n}\bigr) \bigl(1+\mathcal{O}(n^{-1/2})\bigr)
\Bigr).
\end{aligned}
\end{equation*}
Since $|(-1/2)^n| = 2^{-n}$ and $|\exp(\pm\mathrm{i}\sqrt{2}\,a\sqrt{n})| = 1$, we deduce that $|w_n| \leq C 2^{-n} n^{\frac{1}{4}-\frac{\nu}{2}}$ and
\begin{equation*}
|u_n| \leq C \frac{n^{\frac{1}{4}-\frac{\nu}{2}}}{2^n n!}.
\end{equation*}

\end{proof}

\begin{remark}
In view of the connection between $u_n$ and generalized Laguerre polynomials in \eqref{un-Laguerre}, Theorem~\ref{thm-ab} can also be deduced directly from the classical asymptotic formula for Laguerre polynomials \cite[Eq.~(18.15.14)]{Olver-2010}. By contrast, the self-contained derivation presented here directly reveals how the structure of the recurrence governs the asymptotic behavior of $u_n$.
\end{remark}

Let $u_n^*$ denote the exact value defined by the recurrence relation \eqref{eq:un-recurrence} with its initial conditions, and let $\widehat{u}_n$ be the numerical solution computed in floating-point arithmetic.
Since each step introduces a local roundoff error
\begin{equation*}
\delta_n:=\widehat u_n-\frac{a^2-2(\nu+2n-2)}{4n(\nu+n-1)}\widehat u_{n-1}+\frac{1}{4n(\nu+n-1)}\widehat u_{n-2},\qquad n\ge2,
\end{equation*}
a forward error $|\widehat{u}_n - u_n^*|$ is incurred. The following theorem establishes that, for fixed $a$ and $\nu$, $|\widehat{u}_n - u_n^*|$ is uniformly bounded by the local errors up to step $n$.

\begin{theorem}\label{thm-coefficient-error}
Let $a>0$ and $\nu \notin \{0, -1, -2, \dots\}$. Then there exists $C_{a,\nu}>0$ such that
\begin{equation*}
|\widehat u_n-u_n^*|\le C_{a,\nu} \max\{ |\widehat u_0-u_0^*|,|\widehat u_1-u_1^*|,|\delta_2|,\ldots,|\delta_n|\}.
\end{equation*}
\end{theorem}

\begin{proof}
For notational convenience, we define $A_n := (a^2 - 2(\nu + 2n - 2)) / (4n(\nu+n-1))$ and $B_n := 1 / (4n(\nu+n-1))$.
Then the exact solution $u_n^*$ and the computed solution $\widehat{u}_n$ satisfy
\begin{equation*}
u_n^* = A_n u_{n-1}^* - B_n u_{n-2}^* \quad \text{and} \quad \widehat{u}_n = A_n \widehat{u}_{n-1} - B_n \widehat{u}_{n-2} + \delta_n, \qquad n \ge 2.
\end{equation*}
Subtracting the exact recurrence from the computed one and applying the triangle inequality, we obtain
\begin{equation*}
|\widehat{u}_n - u_n^*| \le |A_n| |\widehat{u}_{n-1} - u_{n-1}^*| + |B_n| |\widehat{u}_{n-2} - u_{n-2}^*| + |\delta_n|, \qquad n \ge 2.
\end{equation*}

Define an auxiliary sequence by
\begin{equation*}
H_0 = H_1 = 1, \qquad H_n = |A_n| H_{n-1} + |B_n| H_{n-2} + 1, \quad n \ge 2.
\end{equation*}
We claim that for all $n \ge 0$,
\begin{equation} \label{eq:error_bound_Hn}
|\widehat{u}_n - u_n^*| \le H_n \max\bigl\{ |\widehat{u}_0 - u_0^*|, |\widehat{u}_1 - u_1^*|, |\delta_2|, \dots, |\delta_n| \bigr\}.
\end{equation}
We proceed by induction. For $n=0$ and $n=1$, the inequality holds trivially since $H_0 = H_1 = 1$. Assume that the claim holds for all steps up to $n-1$. At step $n$, we have
\begin{equation*}
\begin{aligned}
|\widehat{u}_n - u_n^*| &\le |A_n| |\widehat{u}_{n-1} - u_{n-1}^*| + |B_n| |\widehat{u}_{n-2} - u_{n-2}^*| + |\delta_n| \\
&\le |A_n| H_{n-1} \max\bigl\{ |\widehat{u}_0 - u_0^*|, |\widehat{u}_1 - u_1^*|, |\delta_2|, \dots, |\delta_{n-1}| \bigr\} \\
&\quad + |B_n| H_{n-2} \max\bigl\{ |\widehat{u}_0 - u_0^*|, |\widehat{u}_1 - u_1^*|, |\delta_2|, \dots, |\delta_{n-2}| \bigr\} \\
&\quad + 1 \cdot \max\bigl\{ |\widehat{u}_0 - u_0^*|, |\widehat{u}_1 - u_1^*|, |\delta_2|, \dots, |\delta_n| \bigr\} \\
&\le \bigl( |A_n| H_{n-1} + |B_n| H_{n-2} + 1 \bigr) \max\bigl\{ |\widehat{u}_0 - u_0^*|, |\widehat{u}_1 - u_1^*|, |\delta_2|, \dots, |\delta_n| \bigr\} \\
&= H_n \max\bigl\{ |\widehat{u}_0 - u_0^*|, |\widehat{u}_1 - u_1^*|, |\delta_2|, \dots, |\delta_n| \bigr\},
\end{aligned}
\end{equation*}
which completes the induction step.

It remains to show that the auxiliary sequence $H_n$ is uniformly bounded. As $n \to \infty$, we have $A_n \sim -1/n \to 0$ and $B_n \sim 1/(4n^2) \to 0$. Consequently, $|A_n| + |B_n| \to 0$. Thus, there exists an integer $m \ge 2$ such that, for all $n \ge m$
\begin{equation*}
|A_n| + |B_n| \le \frac{1}{2}.
\end{equation*}
We define the constant
\begin{equation*}
C_{a,\nu} = \max\bigl\{ 2, H_0, H_1, \dots, H_{m-1} \bigr\}.
\end{equation*}
By definition, $H_n \le C_{a,\nu}$ holds trivially for the initial segment $0 \le n \le m-1$.
Proceeding by mathematical induction for the subsequent steps $n \ge m$, assuming that $H_k \le C_{a,\nu}$ for all $k < n$, we obtain
\begin{equation*}
\begin{aligned}
H_n &= |A_n| H_{n-1} + |B_n| H_{n-2} + 1 \le |A_n| C_{a,\nu} + |B_n| C_{a,\nu} + 1 \\
&= \bigl( |A_n| + |B_n| \bigr) C_{a,\nu} + 1 \le \frac{1}{2} C_{a,\nu} + 1 \le C_{a,\nu},
\end{aligned}
\end{equation*}
where the final inequality holds since $C_{a,\nu} \ge 2$. Substituting $H_n \le C_{a,\nu}$ into the error inequality \eqref{eq:error_bound_Hn} yields the desired result, completing the proof.
\end{proof}

\section{The SE-RR algorithm and error estimates}

In this section, we build on the results from Section~\ref{sec3} to construct the SE-RR algorithm for the generalized Marcum $Q$-function. Furthermore, we analyze its computational complexity, establish a truncation error estimate, and derive an upper bound on the condition number to assess potential cancellation when evaluating the truncated series.

\begin{theorem}
\label{thm:se-rr-series}
Let $a, \nu>0$ and $b\geq0$.
Then
\begin{equation*}
Q_\nu(a,b) = 1 - \frac{e^{-\frac{a^2}{2}}}{2^{\nu-1} \Gamma(\nu)} \sum_{n=0}^\infty u_n(a,\nu) \frac{b^{2 (\nu + n)}}{2 (\nu + n)},
\end{equation*}
where $u_0=1$, $u_1=\frac{a^2 - 2\nu}{4\nu}$ and
\begin{equation*}
u_n = \frac{a^2 - 2(\nu + 2n - 2)}{4n(\nu+n-1)} u_{n-1} - \frac{1}{4n(\nu+n-1)} u_{n-2}, \quad n\geq 2.
\end{equation*}
\end{theorem}

\begin{proof}
By \eqref{Qab-symmetry} and Corollary \ref{cor-main}, we obtain
\begin{equation*}
\begin{aligned}
Q_\nu(a,b) &= 1 - \frac{e^{-\frac{a^2}{2}}}{a^{\nu-1}} \int_{0}^b t^\nu e^{-\frac{t^2}{2}}I_{\nu-1}(at) \textrm{d} t \\
&= 1 - \frac{e^{-\frac{a^2}{2}}}{a^{\nu-1}} \int_{0}^b t^\nu \frac{a^{\nu-1}}{2^{\nu-1} \Gamma(\nu)} \sum_{n=0}^\infty u_n(a,\nu) t^{2n+\nu-1} \textrm{d} t \\
&= 1 - \frac{e^{-\frac{a^2}{2}}}{2^{\nu-1} \Gamma(\nu)} \sum_{n=0}^\infty u_n(a,\nu) \int_{0}^b t^\nu t^{2n+\nu-1} \textrm{d} t \\
&= 1 - \frac{e^{-\frac{a^2}{2}}}{2^{\nu-1} \Gamma(\nu)} \sum_{n=0}^\infty u_n(a,\nu) \frac{b^{2 (\nu + n)}}{2 (\nu + n)}.
\end{aligned}
\end{equation*}
\end{proof}

\begin{remark}
By Remark~\ref{rem:laguerre-coefficients}, with $p=2$ and $\mu=\nu-1$, the coefficients in the preceding theorem satisfy
\begin{equation*}
u_n(a,\nu)=\frac{(-1/2)^n}{(\nu)_n}L_n^{(\nu-1)}\Big(\frac{a^2}{2}\Big).
\end{equation*}
Substitution of this identity into the series in Theorem \eqref{thm:se-rr-series} gives
\begin{equation*}
Q_\nu(a,b)=1-e^{-a^2/2}\sum_{n=0}^{\infty}\frac{(-1)^nL_n^{(\nu-1)}(a^2/2)}{\Gamma(\nu+n+1)}\Big(\frac{b^2}{2}\Big)^{n+\nu},
\end{equation*}
For $a,\nu>0$ and $b\geq0$, this is precisely the orthogonal-polynomial representation derived by Andr\'as, Baricz, and Sun \cite[Eq.~(8)]{Andras-AUSM-2011}. Their work provides an earlier derivation of this series based on a different approach.
\end{remark}

Setting $\nu=1$ yields a formula for computing the Marcum $Q$-function.

\begin{corollary}
Let $a>0$ and $b\geq0$. Then
\begin{equation*}
Q(a,b) = 1 - e^{-\frac{a^2}{2}} \sum_{n=0}^\infty u_n(a) \frac{b^{2n+2}}{2(n+1)},
\end{equation*}
where $u_0=1$, $u_1=\frac{a^2 - 2}{4}$, and the recurrence relation simplifies to
\begin{equation*}
u_n = \frac{a^2 - 4n + 2}{4n^2} u_{n-1} - \frac{1}{4n^2} u_{n-2}, \quad n\geq 2.
\end{equation*}
\end{corollary}

Consequently, we obtain pseudocode for the algorithm computing $Q_\nu(a,b)$, see Algorithm \ref{alg1}.

\begin{algorithm}[htbp]
\caption{Computation of $Q_\nu(a,b)$}
\label{alg1}
\DontPrintSemicolon
\KwIn{$a, b, \nu$; Iterations $N$}
\KwOut{$Q \approx Q_\nu(a,b)$}

\BlankLine
\tcp{1. Precompute coefficients $u_n$}
Initialize array $u[0 \dots N]$\;
$u_0 \leftarrow 1$, \quad $u_1 \leftarrow \frac{a^2 - 2\nu}{4\nu}$\;
\For{$n \leftarrow 2$ \KwTo $N$}{
    $Den \leftarrow 4n(\nu+n-1)$\;
    $u_n \leftarrow \frac{a^2 - 2(\nu + 2n - 2)}{Den} u_{n-1} - \frac{1}{Den} u_{n-2}$\;
}

\BlankLine
\tcp{2. Evaluate Series}
$Sum \leftarrow 0$\;
$PowerTerm \leftarrow b^{2\nu}$

\For{$n \leftarrow 0$ \KwTo $N$}{
    \tcp{Term corresponds to $\frac{b^{2(\nu+n)}}{2(\nu+n)}$}
    $Sum \leftarrow Sum + u_n \cdot \frac{PowerTerm}{2(\nu + n)}$\;
    $PowerTerm \leftarrow PowerTerm \cdot b^2$
}

\BlankLine
\Return{$Q=1 - \frac{e^{-\frac{a^2}{2}}}{2^{\nu-1} \Gamma(\nu)} \cdot Sum$}\;
\end{algorithm}

\vspace{0.5cm}

Next, we briefly analyze the computational complexity of SE-RR. For a truncation order $N$, each step of the coefficient recurrence and series summation requires $\mathcal{O}(1)$ operations. Thus, the total computational complexity scales linearly with the truncation order, denoted as $\mathcal{O}(N)$. Since all coefficients are stored before summation, the space complexity is $\mathcal{O}(N)$.

Now, we analyze the approximation error of $Q_\nu(a,b)$ using the truncated series defined by
\begin{equation*}
Q_\nu^N(a,b) := 1 - \frac{e^{-\frac{a^2}{2}}}{2^{\nu-1}\Gamma(\nu)} \sum_{n=0}^N u_n(a,\nu) \frac{b^{2 (\nu + n)}}{2 (\nu + n)}.
\end{equation*}

\begin{theorem} \label{thm-error}
Let $a>0$, $b\geq0$, $N>0$, and $\nu > 0$. Then, there exists a constant $K_{a,\nu}$ depending on $a$ and $\nu$ such that
\begin{equation*}
| Q_\nu^N(a,b) - Q_\nu(a,b) | \le K_{a,\nu} \frac{e^{\frac{b^2}{2}} b^{2(\nu + N + 1)}}{2^{N} (N+1)! (N+1)^{3/4+\nu/2}}.
\end{equation*}
\end{theorem}

\begin{proof}
Based on the previous analysis, there exists a constant $D_{a,\nu}$ such that $|u_n(a,\nu)| \le D_{a,\nu} \frac{n^{1/4-\nu/2}}{2^n n!}$ for all $n \ge 1$.
The error term is given by the tail of the series:
\begin{equation*}
| Q_\nu^N(a,b) - Q_\nu(a,b) | = \frac{e^{-\frac{a^2}{2}}}{2^{\nu-1}\Gamma(\nu)} \left| \sum_{n=N+1}^\infty u_n(a,\nu) \frac{b^{2 (\nu + n)}}{2 (\nu + n)} \right|.
\end{equation*}
Let $A_{a,\nu} = \frac{e^{-a^2/2}}{2^{\nu-1}\Gamma(\nu)}$. Using the bound for $u_n$ and the triangle inequality, we have
\begin{equation*}
| Q_\nu^N(a,b) - Q_\nu(a,b) | \le A_{a,\nu} \sum_{n=N+1}^\infty D_{a,\nu} \frac{n^{1/4-\nu/2}}{2^n n!} \frac{b^{2 \nu + 2n}}{2(\nu + n)}.
\end{equation*}
For $\nu > 0$ and $n \ge N+1$, we have $n^{1/4-\nu/2}/(\nu+n)\le n^{-3/4-\nu/2}\le (N+1)^{-3/4-\nu/2}$. Thus,
\begin{equation*}
\begin{aligned}
| Q_\nu^N(a,b) - Q_\nu(a,b) | &\le \frac{A_{a,\nu} D_{a,\nu}}{2(N+1)^{3/4+\nu/2}} b^{2\nu} \sum_{n=N+1}^\infty \frac{1}{n!} \Big( \frac{b^2}{2} \Big)^n.
\end{aligned}
\end{equation*}
We apply the bound for the tail of the exponential series: $\sum_{n=N+1}^\infty \frac{x^n}{n!} \le \frac{x^{N+1}}{(N+1)!}e^x$ for $x \ge 0$. Setting $x = b^2/2$, we obtain
\begin{equation*}
\begin{aligned}
| Q_\nu^N(a,b) - Q_\nu(a,b) | &\le \frac{A_{a,\nu} D_{a,\nu}}{2(N+1)^{3/4+\nu/2}} b^{2\nu} \frac{(b^2/2)^{N+1}}{(N+1)!} e^{\frac{b^2}{2}} \\
&= \frac{A_{a,\nu} D_{a,\nu}}{2(N+1)^{3/4+\nu/2}} \frac{b^{2\nu} b^{2N+2}}{2^{N+1} (N+1)!} e^{\frac{b^2}{2}} \\
&= \Big( \frac{A_{a,\nu} D_{a,\nu}}{4} \Big) \frac{e^{\frac{b^2}{2}} b^{2(\nu + N + 1)}}{2^{N} (N+1)! (N+1)^{3/4+\nu/2}}.
\end{aligned}
\end{equation*}
Defining $K_{a,\nu} = \frac{A_{a,\nu} D_{a,\nu}}{4}$, the proof is complete.
\end{proof}

\begin{remark}
We note that Andr\'as, Baricz, and Sun~\cite[Section~2.3]{Andras-AUSM-2011} also studied the truncation error of this series for $Q_\nu(a,b)$ and obtained a similar bound.
\end{remark}

Finally, we investigate cancellation errors in the evaluation of the truncated series. Let $T_n := \frac{e^{-a^2/2}u_n(a,\nu)b^{2(\nu+n)}}{2^\nu\Gamma(\nu)(\nu+n)}$. We approximate the generalized Marcum $Q$-function using the truncated series $Q_\nu^N(a,b)=1 - \sum_{n=0}^{N} T_n$.
Two types of cancellation can happen here. First, the terms $T_n$ can have opposite signs, causing large positive and negative terms to cancel. Second, when $\sum_{n=0}^{N} T_n$ is close to $1$, subtracting it from $1$ can cause severe cancellation, making the final result sensitive to errors in the sum.

To measure this sensitivity, we define the condition number
\begin{equation*}
\kappa = \frac{1 + \sum_{n=0}^{N} |T_n|}{\bigl| 1 - \sum_{n=0}^{N} T_n \bigr|}\geq 1.
\end{equation*}
When strong cancellation occurs, the numerator is much larger than the denominator, making $\kappa$ large. Conversely, when cancellation is weak, the numerator and denominator are of comparable size, so $\kappa \approx 1$. The following theorem provides an explicit upper bound on the condition number $\kappa$, giving a sufficient condition to control the error amplification caused by both internal cancellations and the final subtraction.

\begin{theorem}\label{cancell-bound}
Let $a > 0$, $b \ge 0$, and $\nu > 0$. For all sufficiently large $N$, we have
\begin{equation*}
\kappa \le B_{a,b,\nu} =
\begin{cases}
\dfrac{1 + M_{a,b,\nu}}{1 - M_{a,b,\nu}},
& \text{if } M_{a,b,\nu}<1, \\[8pt]
\dfrac{2\Gamma(\nu)(1+M_{a,b,\nu})}{\Gamma(\nu,b^2/2)},
& \text{if } M_{a,b,\nu}\geq1,
\end{cases}
\end{equation*}
where $M_{a,b,\nu} = \frac{2}{\Gamma(\nu+1)} e^{-a^2/4 + b^2/2} (b^2/2)^\nu$.
\end{theorem}

\begin{proof}
Notice that
\begin{equation*}
T_n = \frac{(-1)^n e^{-a^2/2}}{\Gamma(\nu+n+1)} \Bigl(\frac{b^2}{2}\Bigr)^{\nu+n} L_n^{(\nu-1)}\Bigl(\frac{a^2}{2}\Bigr).
\end{equation*}
Using the classical Laguerre inequality \cite[Eq.~(18.14.8)]{Olver-2010},
\begin{equation*}
\bigl|L_n^{(\nu)}(x)\bigr| \le e^{x/2} \frac{\Gamma(\nu+n+1)}{\Gamma(\nu+1) n!}, \qquad n \ge 0,
\end{equation*}
and applying this bound to the identity $L_n^{(\nu-1)}(x) = L_n^{(\nu)}(x) - L_{n-1}^{(\nu)}(x)$ for $n \ge 1$, we obtain
\begin{equation*}
\Bigl|L_n^{(\nu-1)}\Bigl(\frac{a^2}{2}\Bigr)\Bigr| \le e^{a^2/4} \frac{\Gamma(\nu+n)(\nu+2n)}{\Gamma(\nu+1) n!}.
\end{equation*}
This inequality also holds trivially for $n=0$ since $L_0^{(\nu-1)} \equiv 1$. Then, we have
\begin{equation*}
|T_n| \le \frac{e^{-a^2/4}}{\Gamma(\nu+1)} \frac{\nu+2n}{\nu+n} \frac{(b^2/2)^{\nu+n}}{n!} \le \frac{2 e^{-a^2/4}}{\Gamma(\nu+1)} \frac{(b^2/2)^{\nu+n}}{n!}.
\end{equation*}
Summing this inequality over $n=0, \dots, N$ yields
\begin{equation*}
\sum_{n=0}^{N} |T_n| \le \frac{2 e^{-a^2/4}}{\Gamma(\nu+1)} \Bigl(\frac{b^2}{2}\Bigr)^\nu \sum_{n=0}^{N} \frac{(b^2/2)^n}{n!} \le \frac{2 e^{-a^2/4 + b^2/2}}{\Gamma(\nu+1)} \Bigl(\frac{b^2}{2}\Bigr)^\nu = M_{a,b,\nu}.
\end{equation*}

(I) If $M_{a,b,\nu} < 1$, the reverse triangle inequality ensures that
\begin{equation*}
\Bigl|1 - \sum_{n=0}^{N} T_n \Bigr| \ge 1 - \sum_{n=0}^{N} |T_n| \ge 1 - M_{a,b,\nu} > 0,
\end{equation*}
which implies
\begin{equation*}
\kappa = \frac{1 + \sum_{n=0}^{N} |T_n|}{\bigl| 1 - \sum_{n=0}^{N} T_n \bigr|} \le \frac{1 + M_{a,b,\nu}}{1 - M_{a,b,\nu}}.
\end{equation*}

(II) If $M_{a,b,\nu} \ge 1$, since the sequence $n \mapsto \Gamma(\nu+n, b^2/2)/\Gamma(\nu+n)$ is nondecreasing, it follows from \eqref{QS} that
\begin{equation*}
\begin{aligned}
Q_\nu(a,b) &= e^{-a^2/2} \sum_{n=0}^{\infty} \frac{(a^2/2)^n}{n!} \frac{\Gamma(\nu+n, b^2/2)}{\Gamma(\nu+n)} \\
&\ge \frac{\Gamma(\nu, b^2/2)}{\Gamma(\nu)} e^{-a^2/2} \sum_{n=0}^{\infty} \frac{(a^2/2)^n}{n!} = \frac{\Gamma(\nu, b^2/2)}{\Gamma(\nu)} > 0.
\end{aligned}
\end{equation*}
For all sufficiently large $N$, we have
\begin{equation*}
\bigl| 1 - \sum_{n=0}^{N} T_n \bigr| = \bigl|Q_\nu^N(a,b)\bigr| \ge \frac{1}{2} Q_\nu(a,b) \ge \frac{\Gamma(\nu, b^2/2)}{2\Gamma(\nu)}.
\end{equation*}
Combining this lower bound with $\sum_{n=0}^{N} |T_n| \le M_{a,b,\nu}$ yields the second bound.
\end{proof}

\begin{figure}[H]
    \centering
    \includegraphics[width=\textwidth]{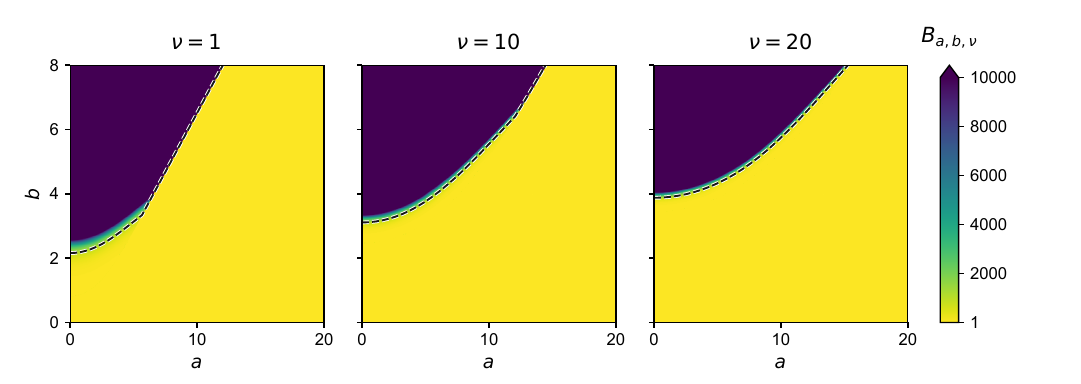}
    \caption{Heatmaps of $B_{a,b,\nu}$ for $\nu=1,10,20$. Dashed curves mark the boundary of $B_{a,b,\nu}\leq10^3$.}
    \label{fig:cancell-bound}
\end{figure}

To illustrate the bound in Theorem~\ref{cancell-bound}, Figure~\ref{fig:cancell-bound} shows heatmaps of $B_{a,b,\nu}$ over the $(a,b)$-plane for $\nu=1,10,20$. Small values of $B_{a,b,\nu}$ ensure that the condition number remains small, as in the lower-right region. In this paper, we regard $B_{a,b,\nu}\leq10^3$ as acceptable for using SE-RR, and use the second formula when $M_{a,b,\nu}\geq1-2\times10^{-3}$.

For example, Table~\ref{tab:cancellation-examples} presents the results for three parameter sets. The floating-point relative error is the relative difference between the double-precision result and the high-precision truncated result, with the latter used as the reference. When $B_{a,b,\nu}$ is small, the actual condition number $\kappa$ and the floating-point relative error are both small. The computation remains accurate even when $B_{a,b,\nu}$ approaches $1000$, since the actual condition number is still small. However, for $(a,b,\nu)=(1,7,1)$, severe cancellation occurs: both $\kappa$ and its upper bound are extremely large, and the double-precision evaluation is unreliable.

\begin{table}[htbp]
    \centering
    \setlength{\tabcolsep}{3pt}
    \renewcommand{\arraystretch}{1.15}
    \resizebox{\textwidth}{!}{
    \begin{tabular}{c c c c c c c l}
        \toprule
        $(a,b,\nu)$ & $B_{a,b,\nu}$ & $\kappa$ & $\sum\lvert T_n\rvert$ & $\sum T_n$ & $1-\sum T_n$ & FP relative error & Remarks \\
        \midrule
        $(8,5,1)$ & $7.16$ & $1.0043$ & $0.003225$ & $0.001040$ & $0.998960$ & $6.88\times10^{-18}$ & Usable \\
        $(3,2.556,1)$ & $999.32$ & $2.1629$ & $0.590251$ & $0.264775$ & $0.735225$ & $2.14\times10^{-17}$ & Usable \\
        $(1,7,1)$ & $1.46\times10^{23}$ & $2.98\times10^{18}$ & $8.01\times10^9$ & $0.9999999973$ & $2.69\times10^{-9}$ & $13.93$ & Not usable \\
        \bottomrule
    \end{tabular}
    }
    \caption{Cancellation and floating-point errors for three parameter sets, with $N=128$.}
    \label{tab:cancellation-examples}
\end{table}

\section{Numerical experiments and results}

\subsection{Numerical experiments and results}

In this section, we present comprehensive numerical experiments to demonstrate the efficiency of the proposed SE-RR algorithm. Specifically, we compare SE-RR against the SE, AF, and NI methods described in Section~2 and implemented in \cite{Gil-ACM-2014}.

We evaluate the efficiency of each method by the computation time required to achieve a target accuracy. To ensure a fair comparison, independent high-precision evaluations based on the definition \eqref{Qab} of the Marcum $Q$-function were used as reference values. We recorded the smallest tested truncation order for SE, SE-RR, and AF, or the smallest tested quadrature refinement level for NI, that achieved a relative error no greater than $\epsilon=10^{-10}$ against this reference. Each method was then rerun with the selected setting to measure its computation time. Reference evaluation, parameter searches, and error checks were excluded from the timing.

To assess performance across different programming environments, we evaluated all four methods separately in Python, C++, and Fortran, with timing comparisons made within each language. All methods used double-precision floating-point arithmetic. In C++ and Fortran, the same optimization settings were applied to all four methods, with performance optimizations enabled and aggressive optimizations that could alter floating-point results disabled. All experiments were conducted on a personal computer equipped with an AMD Ryzen 7 5800H processor and 32\,GB of RAM.

\begin{figure}[H]
    \centering
    \includegraphics[width=0.8\textwidth]{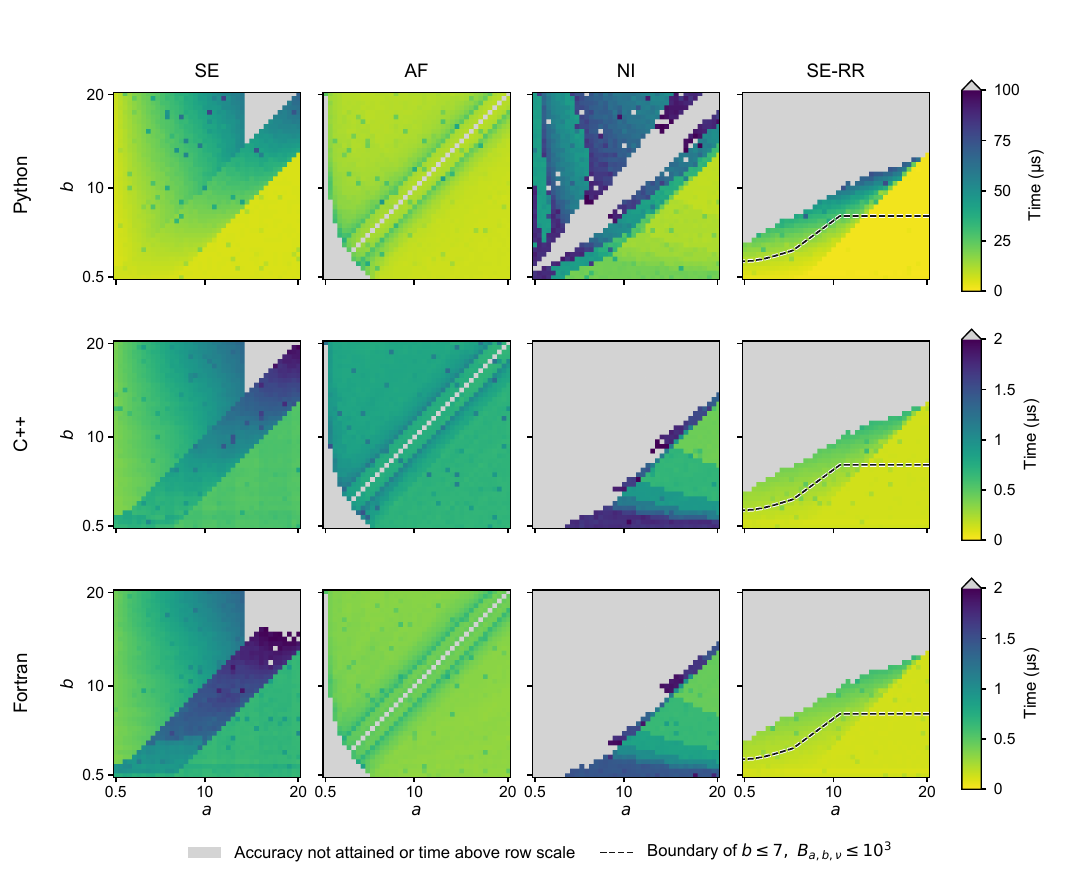}
    \caption{Computation times for $\nu=1$.}
    \label{fig:individual-nu1}
\end{figure}

\begin{figure}[H]
    \centering
    \includegraphics[width=0.8\textwidth]{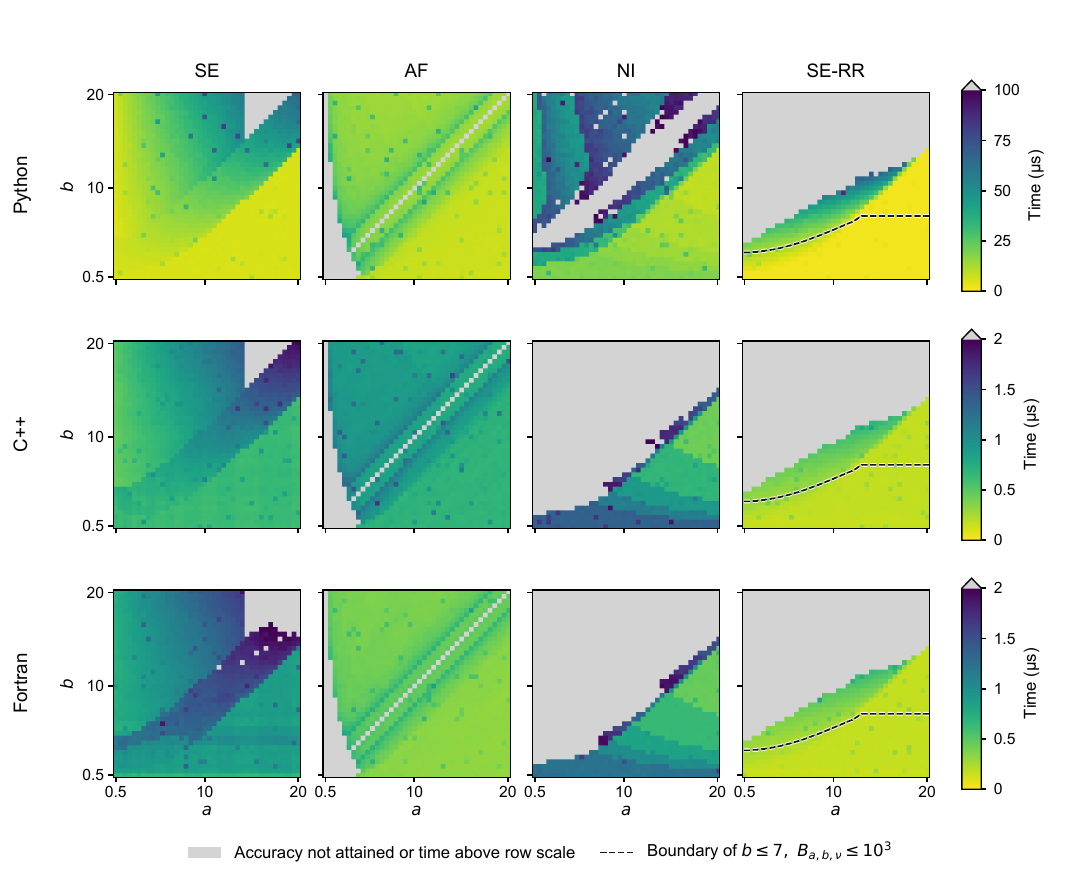}
    \caption{Computation times for $\nu=10$.}
    \label{fig:individual-nu10}
\end{figure}

\begin{figure}[H]
    \centering
    \includegraphics[width=0.8\textwidth]{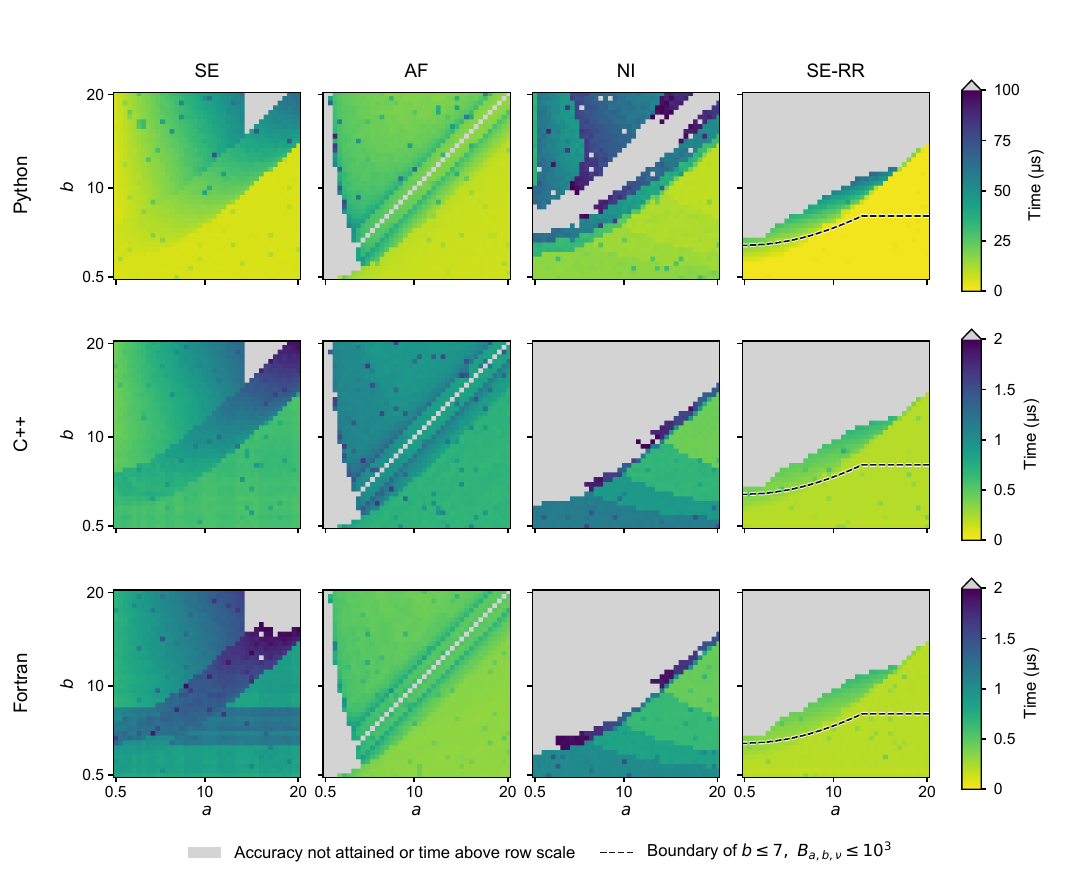}
    \caption{Computation times for $\nu=20$.}
    \label{fig:individual-nu20}
\end{figure}

In each of the three language environments, we tested order $\nu=1,10,20$ on the grid $a,b\in\{0.5,1.0,\ldots,20.0\}$, giving $1{,}600$ points for each $\nu$. Each timing experiment was repeated five times, and the median computation time is reported.
For each language, all methods share a linear color scale capped at $100\,\mu\mathrm{s}$ for Python and $2\,\mu\mathrm{s}$ for C++ and Fortran, where gray points indicate failure to reach the target accuracy or computation times exceeding the upper limit.
Figures~\ref{fig:individual-nu1}--\ref{fig:individual-nu20} display the computation times across this grid. As shown in Figures~\ref{fig:individual-nu1}--\ref{fig:individual-nu20}, SE-RR achieves the fastest computation in the small-$b$ region, particularly for large $a$ and small $b$. However, we also note that SE-RR performs poorly in the large-$b$ region, where it fails to meet the convergence criterion.

To quantify the differences in computational speed, Table~\ref{tab:individual-method-timing} reports the mean computation times of the four methods in each language. We also report the speedup of SE-RR relative to each of the other methods. As shown in the table, SE-RR is faster in most comparisons.

\begin{table}[htbp]
    \centering

    \small
    \setlength{\tabcolsep}{3pt}
    \renewcommand{\arraystretch}{1.15}
    \begin{tabular}{c l r r r r r r r}
        \toprule
        & & \multicolumn{4}{c}{Mean time ($\mu$s)} & \multicolumn{3}{c}{SE-RR speedup} \\
        \cmidrule(lr){3-6}\cmidrule(lr){7-9}
        $\nu$ & Language & SE & AF & NI & SE-RR & vs.\ SE & vs.\ AF & vs.\ NI \\
        \midrule
        1 & Python & 11.3 & \textbf{11.1} & 47 & 12.9 & 0.88 & 0.86 & 3.65 \\
        1 & C++ & 0.691 & 0.756 & 3.08 & \textbf{0.23} & 3.00 & 3.28 & 13.41 \\
        1 & Fortran & 0.927 & 0.386 & 2.83 & \textbf{0.213} & 4.36 & 1.82 & 13.31 \\
        \midrule
        10 & Python & \textbf{9.75} & 12.4 & 35 & 10.2 & 0.95 & 1.21 & 3.42 \\
        10 & C++ & 0.726 & 0.779 & 2.18 & \textbf{0.254} & 2.86 & 3.07 & 8.61 \\
        10 & Fortran & 1.03 & 0.398 & 2.05 & \textbf{0.236} & 4.38 & 1.68 & 8.68 \\
        \midrule
        20 & Python & 8.56 & 13.3 & 20.7 & \textbf{7.86} & 1.09 & 1.69 & 2.64 \\
        20 & C++ & 0.669 & 0.797 & 1.2 & \textbf{0.27} & 2.48 & 2.95 & 4.44 \\
        20 & Fortran & 1.06 & 0.409 & 1.15 & \textbf{0.255} & 4.17 & 1.60 & 4.49 \\
        \bottomrule
    \end{tabular}
    \caption{Mean computation times and speedups on the grid points.}
    \label{tab:individual-method-timing}
\end{table}

\subsection{Improved hybrid algorithm}

Based on Theorem~\ref{cancell-bound} and the numerical results in the preceding subsection, we combine the methods in the following hybrid strategy.
\begin{enumerate}
    \item If $b \leq 7$ and $B_{a,b,\nu} \leq 10^3$, use SE-RR.
    \item If ($b > 7$ and $a \leq 7$) or ($b \leq 7$ and $B_{a,b,\nu} > 10^3$), use SE.
    \item If $ab > 30$ and $\nu^2 < 2ab$, use AF.
    \item If $g_1(a, \nu) < b < g_2(a, \nu)$ and $\nu < 135$, use RR.
    \item If $g_1(a, \nu) < b < g_2(a, \nu)$ and $\nu \ge 135$, use LAE.
    \item Otherwise, use NI.
\end{enumerate}

\begin{figure}[H]
    \centering
    \includegraphics[width=100mm]{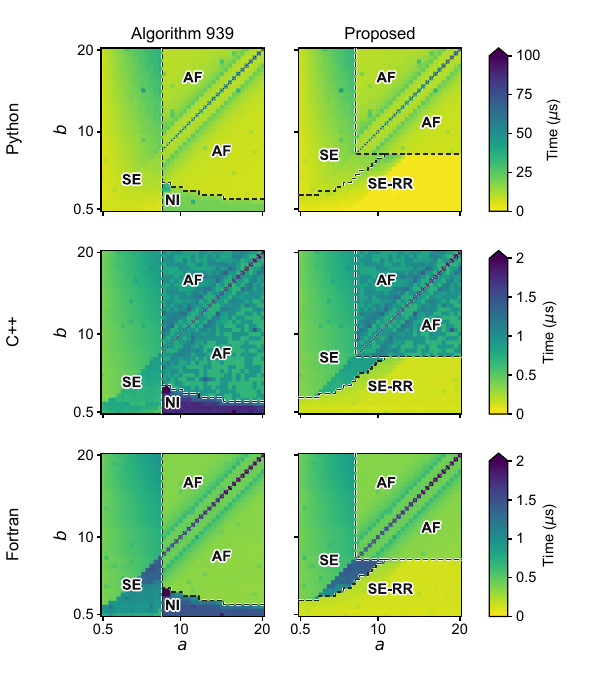}
    \caption{Computation times of Algorithm~939 and the proposed hybrid strategy for $\nu=1$.}
    \label{fig:performance_comparison}
\end{figure}

\begin{figure}[H]
    \centering
    \includegraphics[width=100mm]{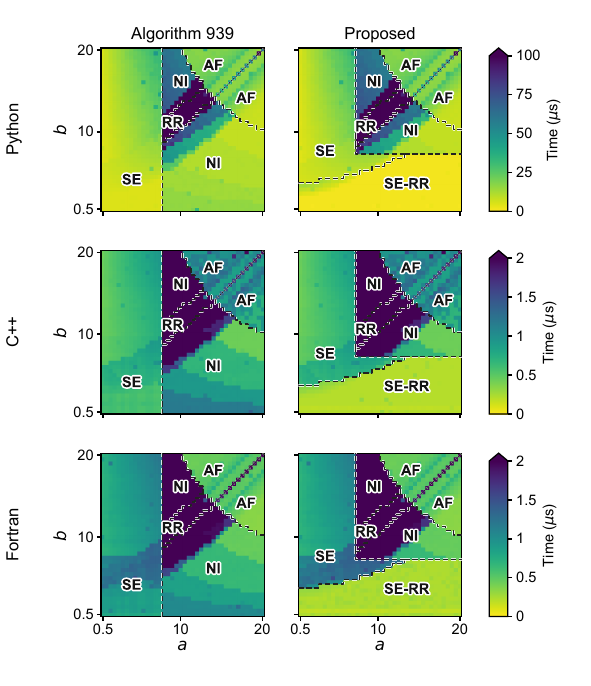}
    \caption{Computation times of Algorithm~939 and the proposed hybrid strategy for $\nu=20$.}
    \label{fig:hybrid-timing-nu20}
\end{figure}

To visually demonstrate the computational advantages of the proposed hybrid strategy, we generate heatmaps illustrating the execution time of $Q_\nu(a,b)$ across the parameter space for fixed $\nu=1,20$ (Figs.~\ref{fig:performance_comparison} and~\ref{fig:hybrid-timing-nu20}) in different programming environments, and compare it with the existing Algorithm~939 proposed in \cite{Gil-ACM-2014}. Specifically, we adopt the officially released Fortran implementation of Algorithm~939~\cite{Gil-ACM-2014-code} as our benchmark and port it directly to Python and C++.

As shown in Figs.~\ref{fig:performance_comparison} and~\ref{fig:hybrid-timing-nu20}, Algorithm~939 uses SE, NI, or AF in the region where $b\leq 7$ and $B_{a,b,\nu}\leq 10^3$. Using SE-RR substantially reduces the computation time in this region, which has the lowest mean computation time among the regions of the proposed hybrid strategy. For example, in Fortran at $\nu=1$, the mean computation time in this region decreases from $0.737$ to $0.155\,\mu\mathrm{s}$, while the mean over all grid points decreases from $0.608$ to $0.458\,\mu\mathrm{s}$.

To further evaluate the performance in Fortran, the most computationally efficient language environment, we generate $10^7$ independent, log-uniformly sampled parameter triples $(a, b, \nu) \in [0.01, 100] \times [0.001, 100] \times [1, 1000]$. The reference values are computed via an independent extended-precision Bessel integral. Subsequently, we perform the evaluations using the open-source implementation of Algorithm~939 and our hybrid algorithm, both targeting a relative accuracy of $10^{-10}$. Overall, the average computation time is $0.929\,\mu\mathrm{s}$ for Algorithm~939 and $0.510\,\mu\mathrm{s}$ for the new hybrid algorithm, yielding an average speedup of $82.2\%$. Specifically, in the region where SE-RR is applied, the average computation time is $0.839\,\mu\mathrm{s}$ for Algorithm~939 and $0.327\,\mu\mathrm{s}$ for the hybrid algorithm, delivering an average speedup of $156.6\%$ and a maximum observed pointwise speedup of $5.49\times$.

\begin{figure}[H]
    \centering
    \includegraphics[width=\linewidth]{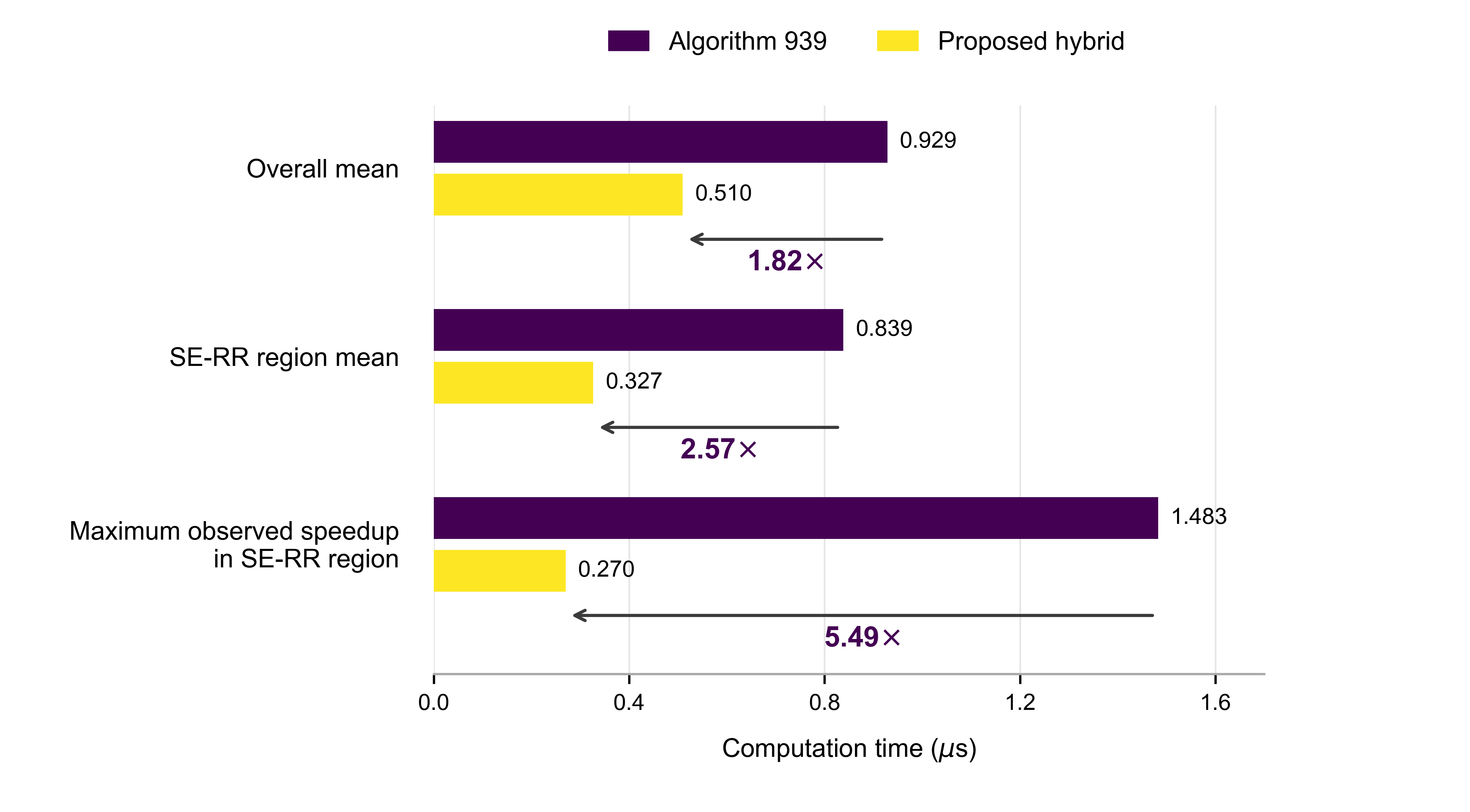}
    \caption{Fortran computation times and speedups of the proposed hybrid algorithm relative to Algorithm~939.}
    \label{fig:random-hybrid-speedup}
\end{figure}

\subsection{Accuracy on random test points}

To test the numerical stability of SE-RR in its intended region of use, where $b\leq 7$ and $B_{a,b,\nu}\leq 10^3$, we conducted a large-scale random test. Specifically, we generated $10^8$ random triples $(a,b,\nu)$ by sampling the three parameters independently and uniformly on a logarithmic scale over $a\in[0.01,100]$, $b\in[0.001,7]$, and $\nu\in[0.1,1000]$. Of these triples, $94{,}671{,}159$ fell within the region $b\leq 7$ and $B_{a,b,\nu}\leq 10^3$. For each of these triples, we used the value computed by the Fortran implementation of Algorithm~939 as the reference value. A triple was counted as passing if SE-RR attained a relative error of at most $10^{-10}$ against this reference value for at least one tested truncation order. Otherwise, it was counted as failing. All $94{,}671{,}159$ test points in this region met the relative error tolerance of $10^{-10}$, giving a pass rate of $100\%$.

\section{Applications}

\subsection{Applications in the computation of noncentral chi distribution}

In this subsection, we demonstrate the practical utility of the proposed algorithm in statistical computing, specifically for evaluating the CDF of the noncentral chi distribution.
Let $R$ denote a random variable following a noncentral chi distribution with $2M$ degrees of freedom and a noncentrality parameter $a$. The CDF of $R$, denoted as $F_R(b)$, is given by
\begin{equation*}
    F_R(b) = \mathbb{P}(R \leq b) = 1 - \mathbb{P}(R > b) = 1 - Q_M(a, b).
\end{equation*}
In this part, we focus on the scenario where $M$ and $a$ are relatively small.

To analyze the statistical properties of $R$ or to visualize its distribution, one frequently needs to evaluate $F_R(b)$ over a dense grid of thresholds $b_1,b_2,\dots,b_N$ with fixed $M$ and $a$.
The proposed hybrid strategy offers two sources of speedup. First, in the parameter configurations considered here, SE-RR replaces SE at thresholds in its region of applicability. Second, for batch evaluation, the coefficients $u_n$ depend only on $M$ and $a$. They are computed once and reused in the series summation at every threshold assigned to SE-RR.

To demonstrate this, we evaluate $F_R(b)$ for six parameter configurations of $M$ and $a$. For each configuration, we use 1000 equally spaced points in $b\in[0,5]$. Gil, Segura, and Temme~\cite{Gil-CPC-2015} built on Algorithm~939 to develop GammaCHI, a Fortran program for computing the noncentral chi-square CDF, and improved the software in 2021~\cite{Gil-CPC-2021}. We use this algorithm as the baseline and measure the time taken by our hybrid algorithm and the baseline to compute the entire curve, see Table~\ref{tab:cdf_time}. Relative to GammaCHI, the proposed hybrid strategy achieves an average speedup of $1.41$ for pointwise evaluation and $1.82$ for batch evaluation.

\begin{table}[htbp]
    \centering

    \small
    \setlength{\tabcolsep}{4pt}
    \begin{tabular}{lrrrrrr}
        \toprule
        & \multicolumn{2}{c}{GammaCHI (ms)} & \multicolumn{2}{c}{Proposed hybrid (ms)} & \multicolumn{2}{c}{Speedup} \\
        \cmidrule(lr){2-3}\cmidrule(lr){4-5}\cmidrule(lr){6-7}
        $(M,a)$ & Pointwise & Batch & Pointwise & Batch & Pointwise & Batch \\
        \midrule
        $(1.0,0.5)$ & 0.472 & 0.352 & 0.340 & 0.194 & 1.39 & 1.81 \\
        $(1.0,2.0)$ & 0.607 & 0.459 & 0.380 & 0.221 & 1.60 & 2.08 \\
        $(3.0,1.0)$ & 0.693 & 0.509 & 0.534 & 0.340 & 1.30 & 1.50 \\
        $(1.5,3.0)$ & 0.761 & 0.580 & 0.543 & 0.321 & 1.40 & 1.81 \\
        $(2.0,5.0)$ & 0.870 & 0.672 & 0.625 & 0.347 & 1.39 & 1.94 \\
        $(3.0,4.0)$ & 0.857 & 0.659 & 0.627 & 0.368 & 1.37 & 1.79 \\
        \bottomrule
    \end{tabular}
    \caption{Fortran runtimes for evaluating 1000 CDF values per parameter configuration.}
    \label{tab:cdf_time}
\end{table}

\subsection{Application in the computation of outage probability of URLLC}
In millimeter-wave communication scenarios with a dominant LoS component, the small-scale fading is commonly modeled by the Rician distribution with factor $K$.
A key reliability metric in URLLC systems is the outage probability, defined as the probability that the instantaneous signal-to-noise ratio (SNR) $\gamma$ falls below a threshold $\gamma_{\text{th}}$. For a normalized average SNR $\bar{\gamma}=1$, the outage probability is given by
\begin{equation*}
P_{\text{out}}(\gamma_{\text{th}})
=1-Q_1\!\Big(\sqrt{2K},
\sqrt{2(1+K)\gamma_{\text{th}}}\Big),
\end{equation*}
where $K$ denotes the Rician factor.

In scenarios dominated by strong LoS components, the non-centrality parameter $a=\sqrt{2K}$ becomes significantly large. Conversely, to maintain the stringent reliability required by URLLC, the outage threshold $\gamma_{\text{th}}$ is typically set low relative to the average SNR, yielding a comparatively small argument $b=\sqrt{2(1+K)\gamma_{\text{th}}}$. This specific regime causes the traditional series expansion method SE to converge slowly or suffer from severe numerical instability. Furthermore, under these conditions, the value of $Q_1\big( a, b \big)$ approaches unity. Consequently, direct computation of the outage probability via $1 - Q_1$ leads to catastrophic cancellation errors in floating-point arithmetic due to the loss of significant digits.

By employing our proposed method, SE-RR, which inherently evaluates the complementary series, $P_{\text{out}}$ can be computed directly and accurately as:
\begin{equation*}
    P_{\text{out}} \approx e^{-K} \sum_{n=0}^N v_n(K) \frac{\big( 2(1+K)\gamma_{\text{th}} \big)^{n+1}}{2 (n+1)},
\end{equation*}
where $v_0=1$, $v_1=\frac{K - 1}{2}$, and the coefficients satisfy the recurrence relation
\begin{equation*}
v_n = \frac{K - 2n + 1}{2n^2} v_{n-1} - \frac{1}{4n^2} v_{n-2}, \quad n \geq 2.
\end{equation*}
This formulation not only robustly handles scenarios with extremely large $K$ factors but also completely circumvents the subtraction of nearly identical values. Thus, our method provides a highly efficient, real-time algorithm for URLLC link adaptation and performance analysis.

To illustrate this, we consider numerical examples with Rician factors $K\in\{5,6,\ldots,35\}\cup\{29.5,30.5\}$ and thresholds $\gamma_{\text{th}}\in\{0.1,0.01,0.001\}$. All 99 parameter pairs lie in the region where the proposed hybrid algorithm uses SE-RR. Accordingly, we use Algorithm~939 as the baseline. The mean computation time per point was $1.649\,\mu\mathrm{s}$ for Algorithm~939 and $0.174\,\mu\mathrm{s}$ for our hybrid algorithm, corresponding to a speedup of $9.48$.

\section{Conclusions}

Driven by the extensive applications of the generalized Marcum $Q$-function, $Q_\nu(a,b)$, its efficient computation remains a focal point of research in applied mathematics and engineering. Conventionally, series expansion methods are favored when the parameter $a$ is relatively small. However, in the regime where the non-centrality parameter $a$ is large and the argument $b$ is small, these traditional expansions suffer from slow convergence and numerical instability. Consequently, state-of-the-art implementations are often forced to rely on computationally intensive numerical integration to maintain accuracy.

Motivated by this limitation, we developed the SE-RR algorithm to specifically target this challenging regime. By constructing a Maclaurin series with respect to $b$ rather than $a$, our approach efficiently evaluates cases characterized by a large non-centrality parameter $a$ and a small argument $b$. The introduction of SE-RR successfully overcomes the critical computational bottleneck inherent in existing gold-standard Algorithm 939, delivering a $2.57\times$ speedup in the ``large $a$, small $b$'' regime. By integrating this highly efficient technique, we have established a comprehensive, robust, and fast composite algorithm capable of evaluating the generalized Marcum $Q$-function seamlessly across the entire parameter space.

A direction for future work is to investigate whether this approach can be extended to other integrals involving special functions whose integrands admit series expansions with recursively computable coefficients.

\section*{Acknowledgments}
The authors sincerely thank the editor and the anonymous reviewers for their valuable comments and constructive suggestions, which have greatly improved the clarity, rigor, and presentation of this paper.

\end{document}